\documentclass[12pt,leqno]{amsart}
\usepackage{latexsym,amsmath,amssymb}
\usepackage{graphicx}

\title[A characterization of BLD maps]{A new characterization of the mappings of bounded length distortion}

\author{Piotr Haj{\l}asz, Soheil Malekzadeh}

\address{P.\ Haj{\l}asz: Department of Mathematics, University of Pittsburgh, 301
  Thackeray Hall, Pittsburgh, PA 15260, USA, {\tt hajlasz@pitt.edu}}

\address{S. Malekzadeh: Department of Mathematics, University of Pittsburgh, 301
  Thackeray Hall, Pittsburgh, PA 15260, USA, {\tt som13@pitt.edu}}

\thanks{P.H.\ was supported by NSF grant DMS-1161425.}

plus 6pt minus 12pt
\def\rank{{\rm rank\,}}

\newcommand{\cone}{{\text{Cone}}\, }

\def\H{{\mathcal H}}

\newtheorem{theorem}{Theorem}

\theoremstyle{definition}
\newtheorem{remark}[theorem]{Remark}

\newcommand{\barint}{
\rule[.036in]{.12in}{.009in}\kern-.16in \displaystyle\int }

\newcommand{\barcal}{\mbox{$ \rule[.036in]{.11in}{.007in}\kern-.128in\int $}}

\newcommand{\bbbn}{\mathbb N}
\newcommand{\bbbr}{\mathbb R}

\newcommand{\Heis}{\mathbb H}

\def\mvint_#1{\mathchoice
          {\mathop{\vrule width 6pt height 3 pt depth -2.5pt
                  \kern -8pt \intop}\nolimits_{\kern -3pt #1}}%
          {\mathop{\vrule width 5pt height 3 pt depth -2.6pt
                  \kern -6pt \intop}\nolimits_{#1}}%
          {\mathop{\vrule width 5pt height 3 pt depth -2.6pt
                  \kern -6pt \intop}\nolimits_{#1}}%
          {\mathop{\vrule width 5pt height 3 pt depth -2.6pt
                  \kern -6pt \intop}\nolimits_{#1}}}

\numberwithin{theorem}{section} \numberwithin{equation}{section}

\begin{document}


\begin{abstract}
In this paper, we present a new characterization of the mappings of bounded length distortion (BLD for short).
In the original geometric definition it is assumed that a BLD mapping is open, discrete and sense preserving. We prove
that the first two of the three conditions are redundant and the sense-preserving condition can be replaced by
a weaker assumption that the Jacobian is non-negative.
\end{abstract}

\maketitle

\section{Introduction}
\label{introduction}

The class of mappings of bounded length distortion (BLD for short) was introduced by Martio and V\"ais\"al\"a \cite{MV}, and it plays a fundamental role in the 
contemporary development of geometric analysis and geometric topology, especially in the context of branched coverings of metric spaces, see e.g. 
\cite{drasinp,heinonenk,HK,HKM,heinonenr,heinonenr2,heinonens,Ledonne,ledonnep,pankka}.

\noindent
{\bf Analytic definition.}
A mapping $f:\Omega\to\bbbr^n$ defined on a domain $\Omega\subset\bbbr^n$ is said to be of the 
{\em $M$-bounded length distortion ($M$-BLD)} if it is locally Lipschitz, has non-negative Jacobian $J_f\geq 0$ and
\begin{equation}
\label{lr}
M^{-1}\vert h\vert\leq\vert Df(x)h\vert \leq M\vert h\vert,
\quad
\text{for a.e. $x\in\Omega$ and for all $h\in\bbbr^n$.}
\end{equation}

Martio and V\"ais\"al\"a proved that this definition is equivalent to a more geometric one.

\noindent
{\bf Geometric definition.}
A continuous map $f:\Omega\to\bbbr^n$ defined on a domain $\Omega\subset\bbbr^n$ is $M$-BLD
if it is open, discrete, sense-preserving and for any curve $\gamma$ in $\Omega$ we have
\begin{equation}
\label{lr2}
M^{-1}\ell(\gamma)\leq\ell(f\circ\gamma)\leq M\ell(\gamma),
\end{equation}
where $\ell(\gamma)$ denotes the length of the curve $\gamma$.

A mapping is said to be BLD if it is $M$-BLD for some $M\geq 1$.

Let us recall some topological terminology here. A mapping $f:\Omega\to\bbbr^n$ is called {\em open} if it maps open subsets of $\Omega$ onto open subsets 
of $\bbbr^n$. It is {\em discrete} if the inverse image of any point in $\bbbr^n$ is a discrete set in $\Omega$. We say that $f$ is {\em sense-preserving} 
({\em weakly sense-preserving}) if it is 
continuous and the topological degree with respect to any subdomain $D\Subset\Omega$ satisfies: $\deg(f, D, y) > 0$ 
($\deg(f, D, y) \geq 0$) for all $y\in f(D)\setminus f(\partial D)$. Fore more details, see e.g. \cite{rickman}.

The proof of the equivalence of both definitions given in \cite{MV} goes as follows.

If $f$ satisfies the analytic definition, then $f$ is locally $M$-Lipschitz and the right inequality at \eqref{lr2} follows.
Also $f$ is a non-constant quasiregular mapping, and hence by Reshetnyak's theorem \cite{reshetnyak2},
it is open, discrete and sense-preserving. Then Martio and V\"ais\"al\"a concluded the left inequality at \eqref{lr2} using the so called path-lifting argument
which is available for mappings that are open and discrete.

Now suppose that $f$ satisfies the geometric definition. It easily follows from \eqref{lr2} that $f$ is locally $M$-Lipschitz and hence $|Df(x)h|\leq M|h|$
for almost all $x\in\Omega$ and all $h\in\bbbr^n$. Since $f$ is locally Lipschitz, open, discrete and sense-preserving it follows from 
\cite[I.4.11]{rickman} that $J_f\geq 0$. It remains to prove that $|Df(x)h|\geq M^{-1}|h|$. Again, Martio and V\"ais\"al\"a employed the path-lifting 
argument in the proof. This argument could be used only because it was assumed that the mapping was open and discrete.

The above sketch of the proof shows that the path-lifting argument was used
in both directions of the proof of the equivalence of the analytic and the geometric definitions. In words of Martio and V\"ais\"al\"a
{\em (path-lifting) is perhaps the most important tool in the theory of BLD maps.} 

The main observation in our paper is that the implication from \eqref{lr2} to \eqref{lr} does not require the path-lifting argument.
\begin{theorem}
\label{argument}
If a continuous map $f:\Omega\to\bbbr^n$ defined on a domain $\Omega\subset\bbbr^n$ is such that for any 
rectifiable curve $\gamma$ in $\Omega$ inequality \eqref{lr2} is satisfied, then $f$ is locally Lipschitz and \eqref{lr}
is true.
\end{theorem}
The mappings to which Theorem~\ref{argument} applies may change orientation. For example it applies to
the folding of the plane $f(x,y)=(|x|,y)$. This map preserves lengths of all curves, yet it changes orientation and hence it is not BLD.
In order to apply Theorem~\ref{argument} to the class of BLD mapping we need a condition that would eliminate mappings like folding of the plane;
it suffices to assume that the Jacobian is non-negative.
The following result which is the main result in our paper is a straightforward consequence of Theorem~\ref{argument}.
\begin{theorem}
\label{main}
A continuous mapping $f:\Omega\to\bbbr^n$ is $M$-BLD, if and only if 
$$
M^{-1}\ell(\gamma)\leq \ell(f\circ\gamma)\leq M\ell(\gamma),
$$
for all rectifiable curves $\gamma$ in $\Omega$ and $Jf\geq 0$ a.e. in $\Omega$.
\end{theorem}
\begin{remark}
The condition about non-negative Jacobain is very natural. Indeed, 
if $f:\Omega\to\bbbr^n$ is a Lipschitz map, that is open, discrete, and sense-preserving, then the Jacobian $Jf\geq 0$ is non-negative a.e. \cite[Lemma~I.4.11]{rickman}.
On the other hand, a Lipschitz mapping $f:\Omega\to\bbbr^n$ with $Jf\geq 0$ a.e. is weakly sense-preserving, \cite[Lemma~VI.5.1]{rickman}. 
\end{remark}
\begin{remark}
Theorem~\ref{main} provides a new and a simpler version of the geometric definition by showing that the strong topological assumptions about
openness and discreteness in the geometric definition of the BLD mappings are redundant.
We want to emphasize that the proof that the analytic definition implies our new geometric definition involves Reshetnyak's theorem and 
the path-lifting argument as described above. Only the other implication is based on a new argument that avoids topological assumptions.
\end{remark}
\begin{remark}
The geometric definition of Martio and V\"ais\"al\"a is the base of the independent theory of BLD mappings between generalized metric manifolds, but
our result applies to the Euclidean setting only.
\end{remark}

Our proof is short and elementary, but we arrived to this simple argument through a rather complicated way by studying unrectifiability of the Heisenberg groups.

A mapping $f:X\to Y$ between metric spaces is said to be a {\em weak BLD} mapping, if there is a constant $M\geq 1$ such that for all {\em rectifiable} curves $\gamma$ in $X$, we have 
\begin{equation}
\label{wBLD}
M^{-1}\ell_X(\gamma)\leq\ell_Y(f\circ\gamma)\leq M\ell_X(\gamma).
\end{equation}
This definition was introduced in \cite{HMz}. See also \cite{Ledonne}, where a stronger condition that \eqref{wBLD} is true for {\em all} curves $\gamma$ is required.

It is a well-known fact that the identity map ${\rm id}\,:\Heis^n\to\bbbr^{2n+1}$ has (locally) the weak BLD property. Here, $\Heis^n$ is the standard Heisenberg group.
One of the results in \cite{HMz} provides a characterization of pure $k$-unrectifiability of a metric space $X$ under the assumption that there is a weak BLD 
mapping $\Phi:X\to\bbbr^N$. Since the identity mapping $\Phi={\rm id}\, :\Heis^n\to\bbbr^{2n+1}$ has the weak BLD property, the result provides a new proof of pure 
$k$-unrectifiability of $\Heis^n$ when $k>n$. It also follows from this characterization of pure $k$-unrectifiability of $X$ that if $f:\Omega\subset\bbbr^n\to\bbbr^m$ 
has the weak BLD property, then $m\geq n$ and $\rank Df=n$ a.e. In particular, weak BLD mappings $f:\Omega\subset\bbbr^n\to\bbbr^n$ satisfy $\vert Jf\vert>0$ a.e. 
We quickly realized that a stronger quantitative estimate $|Jf|\geq C>0$ would imply that weak BLD mappings with $Jf\geq 0$ are BLD which is our Theorem~\ref{main}.
We could prove this estimate using the methods of \cite{HMz}, but the proof was long and complicated; eventually we discovered a very simple argument which is presented in this paper.

Our notation is pretty standard. By $\H^{k}$ we will denote the $k$-dimensional Hausdorff measure. We will also use $\H^n$ to denote the Lebesgue measure in $\bbbr^n$.
By $S^{n-1}$ we will denote the unit sphere in $\bbbr^n$ and $\chi_E$ will stand for the characteristic function of a set $E$.

\noindent
{\bf Acknowledgements.} We would like to thank the referee whose comments and suggestions led to an
improvement of the presentation of the paper.

\section{Proof of Theorem~\ref{argument}}
\label{theproof}
Suppose that $f$ satisfies the assumptions of Theorem~\ref{argument}.

For any $y\in B(x,r)\subset\Omega$, the segment $\overline{xy}$ is mapped onto a curve of length bounded by $M\vert x-y\vert$. 
Hence $\vert f(x) - f(y)\vert\leq M\vert x - y\vert$ which means $f$ is locally $M$-Lipschitz in $\Omega$. So, it is 
differentiable a.e. by Rademacher's theorem and 
$$
\vert Df(x)h \vert \leq M\vert h \vert
\quad
\text{for a.e. $x\in\Omega$ and all $h\in\bbbr^n$.}
$$
It remains to prove that $M^{-1}\vert h \vert \leq \vert Df(x)h \vert$ for a.e. $x\in\Omega$ and all $h\in\bbbr^n$
which is equivalent to proving that
\begin{equation}
\label{lower}
\vert Df(x)h \vert \geq M^{-1}
\end{equation}
for a.e. $x\in\Omega$ and all $h\in S^{n-1}$.

Since $Df$ is a measurable mapping, there exists, for any $m\in\bbbn$, a closed set $K_m\subset\Omega$ such that 
$f$ is differentiable at all points of $K_m$,
$\H^n(\Omega\setminus K_m)<m^{-1}$ and $Df$ is continuous on $K_m$. Since we can exhaust $\Omega$ with the sets $K_m$ up to a set of measure zero, 
it suffices to prove \eqref{lower} for any given $m\in\bbbn$, almost every $x\in K_m$ and all $h\in S^{n-1}$.
Fix $m\in\bbbn$. It suffices to prove inequality \eqref{lower} when
$x\in K_m$ is a density point of $K_m$ and $h\in S^{n-1}$ is arbitrary.
To the contrary suppose that there is a density point $x\in K_m$ and
$h_0\in S^{n-1}$ such that
$$
\vert Df(x)h_0\vert = \alpha<M^{-1}.
$$
Without loss of generality, we may assume that $x=0$. This is only used to simplify notation.

It suffices to show that there is a rectifiable curve $\gamma$ passing through $x$ such that $\ell(f\circ\gamma)<M^{-1}\ell(\gamma)$.
This will contradict \eqref{lr2}. To this end it suffices to show that there are
\begin{equation}
\label{Asia1}
\text{$0<\alpha_1<\alpha_2<M^{-1}$, $v\in S^{n-1}$ and $R>0$}
\end{equation}
such that for the curve
$$
\gamma:[0,R]\to\Omega,
\quad
\gamma(t)=x+tv=tv
$$
the following conditions are satisfied
\begin{equation}
\label{Asia2}
|Df(tv)v|<\alpha_1
\quad
\text{whenever $tv\in K_m$ and $0\leq t\leq R$,}
\end{equation}
\begin{equation}
\label{Asia3}
\H^1(\{t\in [0,R]:\, tv\not\in K_m\})<R\, \frac{\alpha_2-\alpha_1}{M}.
\end{equation}
Indeed, since $f$ is locally $M$-Lipschitz and $\gamma$ is $1$-Lipschitz, $|(f\circ\gamma)'(t)|\leq M$ for a.e. $t\in [0,R]$. Moreover
\eqref{Asia2} implies that
$$
|(f\circ\gamma)'(t)|=|Df(\gamma(t))\gamma'(t)|<\alpha_1
\quad
\text{whenever $\gamma(t)\in K_m$.}
$$
Hence
\begin{eqnarray*}
\ell(f\circ\gamma)
& = &
\int_0^R|(f\circ\gamma)'(t)|\, dt 
 < 
\int_{\{\gamma\not\in K_m\}} M \, dt+
\int_{\{\gamma\in K_m\}}{\alpha_1}\, dt\\ 
& < & 
MR\frac{\alpha_2-\alpha_1}{M} + \alpha_1 R
 = 
R\alpha_2 < RM^{-1} = M^{-1}\ell(\gamma).
\end{eqnarray*}
Therefore it remains to prove that there are $\alpha_1,\alpha_2,v$ and $R$ as in \eqref{Asia1} for which the conditions
\eqref{Asia2} and \eqref{Asia3} are satisfied. 
 
Fix any numbers $\beta,\alpha_1,\alpha_2$ such that $\alpha<\beta<\alpha_1<\alpha_2<M^{-1}$. Observe that $M\geq 1$.
 
Let $\cone(r,\delta)$ be the set of all vectors $h\in \overline{B}(0,r)$ such that the angle between $h$ and $h_0$ is less than or equal $\delta$.

Since $Df(0):\bbbr^n\to\bbbr^n$ is linear and continuous, there exists $\delta>0$ such that
$$
\vert Df(0)h\vert<\beta
\quad
\text{for all $h\in\cone(1,\delta)$.}
$$
Since $Df|_{K_m}$ is continuous at $x=0$, there exists $\tau>0$ such that 
\begin{equation}
\label{Michal}
\vert Df(y)h\vert<\alpha_1,
\quad
\text{for $y\in B(0,\tau)\cap K_m$ and $h\in\cone(1,\delta)$.}
\end{equation}
Let $C'(n,\delta)=\H^n(\cone(1,\delta))$.
By a scaling argument
$$
\H^n(\cone(r,\delta)) = \H^n(\cone(1,\delta))r^n=C'(n,\delta)r^n.
$$
Since $x=0$ is a density point of $K_m$ and $0<1-(\alpha_2-\alpha_1)^n/M^n<1$,
there is $0<R<\tau$ such that $B(0,R)=B(x,R)\subset\Omega$ and 
\begin{equation}
\label{cap}
\H^n(\cone(R,\delta)\cap K_m)>C'(n,\delta)R^n\left(1-\frac{(\alpha_2-\alpha_1)^n}{M^n}\right).
\end{equation}
Now we claim that there is a vector $v=\bar{h}/|\bar{h}|\in S^{n-1}$ for some
$0\neq \bar{h}\in\cone(R,\delta)$ for which \eqref{Asia3} is satisfied.
Indeed, suppose to the contrary that for every $0\neq\bar{h}\in\cone(R,\delta)$ we have
$$
\H^1\big(\big\{ t\in [0,R]:\, t\bar{h}/|\bar{h}|\not\in K_m\big\}\big) \geq R\, \frac{\alpha_2-\alpha_1}{M}\, ,
$$
i.e.
\begin{equation}
\label{curve}
\H^1\big(\big\{ t\in [0,R]:\, t\bar{h}/|\bar{h}|\in K_m\big\}\big) \leq R\, \left(1-\frac{\alpha_2-\alpha_1}{M}\right)\, .
\end{equation}
Hence for any $0\neq\bar{h}\in\cone(R,\delta)$ and $u=\bar{h}/|\bar{h}|\in S^{n-1}$
$$
\int_0^R\chi_{K_m}(tu) t^{n-1}\, dt \leq
\int_{R(\alpha_2-\alpha_1)/M}^R t^{n-1}\, dt.
$$
Thus integration in spherical coordinates yields
\begin{eqnarray*}
\H^n(\cone(R,\delta)\cap K_m)  
& = &
\int_{S^{n-1}\cap\cone(1,\delta)} \int_0^R \chi_{K_m}(tu) t^{n-1}\, dt\,d\H^{n-1}(u) \\
& \leq &
\int_{S^{n-1}\cap\cone(1,\delta)}
\int_{R(\alpha_2-\alpha_1)/M}^R t^{n-1}\, dt\, d\H^{n-1}(u) \\
& = &
\H^n\left(\cone(R,\delta)\setminus\cone\left(R\, \frac{\alpha_2-\alpha_1}{M},\delta\right)\right) \\
& = &
C'(n,\delta)R^n\left(1-\frac{(\alpha_2-\alpha_1)^n}{M^n}\right)
\end{eqnarray*}
which clearly contradicts \eqref{cap}. 

We proved that 
for some $0\neq \bar{h}\in\cone(R,\delta)$
the vector $v=\bar{h}/|\bar{h}|\in S^{n-1}$ satisfies \eqref{Asia3}.
Now \eqref{Asia2} easily follows from \eqref{Michal} with $y=tv$ and $h=v$. Indeed, $v\in\cone(1,\delta)$ and 
$tv\in B(0,R)\subset B(0,\tau)$.
The proof is complete. This also completes the proof of Theorem~\ref{main}.

\end{document}